\documentclass{amsart}

\usepackage{amsthm}
\usepackage{amssymb}
\usepackage[all]{xy}
\usepackage{hyperref}
\usepackage{verbatim}

\SelectTips{cm}{}
\calclayout

\newtheorem{theorem}{Theorem}[section]
\newtheorem{lemma}[theorem]{Lemma}
\newtheorem{proposition}[theorem]{Proposition}
\newtheorem{corollary}[theorem]{Corollary}

\theoremstyle{definition}
\newtheorem{definition}[theorem]{Definition}
\theoremstyle{remark}

\newtheorem{example}[theorem]{Example}

\newtheorem{acknowledgement}{Acknowledgement}

\newcommand{\Ass}{\operatorname{Ass}}

\newcommand{\Spec}{\operatorname{Spec}}

\newcommand{\Ht}{\operatorname{ht}}

\newcommand{\V}{\operatorname{V}}

\newcommand{\Supp}{\operatorname{Supp}}

\newcommand{\Ann}{\operatorname{Ann}}

\newcommand{\Max}{\operatorname{Max}}

\newcommand{\m}{\frak{m}}
\newcommand{\n}{\frak{n}}
\newcommand{\p}{\frak{p}}
\newcommand{\q}{\frak{q}}
\newcommand\ab{\operatorname{\frak b}}

\begin{document}

\author[K. Bahmanpour and K. Divaani-Aazar]{Kamal Bahmanpour and Kamran Divaani-Aazar}

\title[Weakly Laskerian rings versus Noetherian rings]
{Weakly Laskerian rings versus Noetherian rings}

\address{K. Bahmanpour, Faculty of Mathematical Sciences, Department of Mathematics, University
of Mohaghegh Ardabili, 56199-11367, Ardabil, Iran-and-School of Mathematics, Institute for Research
in Fundamental Sciences (IPM), P.O. Box 19395-5746, Tehran, Iran}
\email{bahmanpour.k@gmail.com}

\address{K. Divaani-Aazar, Department of Mathematics, Alzahra University, Vanak, Post Code 19834,
Tehran, Iran-and-School of Mathematics, Institute for Research in Fundamental Sciences (IPM), P.O.
Box 19395-5746, Tehran, Iran}
\email{kdivaani@ipm.ir}


\keywords {FSF modules; ring extensions; trivial extensions; weakly Laskerian modules.\\
The research of the authors are supported by grants from IPM (no. 93130022 and no. 93130212; respectively).}

\begin{abstract} Let $R$ be a commutative ring with identity.  We investigate some ring-theoretic properties
of weakly Laskerian $R$-modules. Our results indicate that weakly Laskerian rings behave as Noetherian ones
in many respects. However, we provide some examples to illustrate the strange behavior of these rings in some
other respects.
\end{abstract}

\maketitle

\tableofcontents

\section{Introduction}

Throughout this paper, all rings are assumed to be commutative with identity. Also, all modules are assumed to be left unitary.

Let $R$ be a ring. An $R$-module $M$ is said to be Laskerian if the zero submodule of every quotient of $M$ has a primary
decomposition. Clearly, any Noetherian $R$-module is Laskerian. As a generalization of this notion, the notion of weakly Laskerian
modules was introduced by the present second author and Mafi in \cite{DM1}. An $R$-module $M$ is said to be weakly Laskerian if
every quotient module of $M$ has finitely many associated prime ideals. The class of weakly Laskerian $R$-modules obviously includes
all Laskerian modules. In Example \ref{Q2}, we provide an example of a non-Laskerian ring which is weakly Laskerian. The class of
weakly Laskerian $R$-modules is large enough to contain all Noetherian and Artinian $R$-modules. One may easily check that it
is a Serre class. This means that in any short exact sequence of $R$-modules and $R$-homomorphisms, the middle module is weakly
Laskerian if and only if the two other modules are weakly Laskerian. In the case $R$ is Noetherian, the present first author proved
that an $R$-module $M$ is weakly Laskerian if and only if it is FSF; see \cite[Theorem 3.3]{Ba}.  Recall that by Quy's definition
\cite[Definition 2.1]{Q}, an $R$-module $M$ is said to be FSF if it possesses a finitely generated submodule $N$ such that
$\Supp_RM/N$ is a finite set.

Let us for a while assume that $R$ is Noetherian. The study of finiteness properties of local cohomology modules of finitely
generated $R$-modules has been an active area of research in recent years. Although, the class of weakly Laskerian $R$-modules
is much larger than that of finitely generated $R$-modules, the analogues of many nice finiteness properties of local cohomology
modules of finitely generated $R$-modules have been established for weakly Laskerian $R$-modules. So, this class deserves
a deeper investigation. In fact, in several papers the class of weakly Laskerian $R$-modules have been examined in conjunction
with local cohomology modules; see e.g. \cite{DM1}, \cite{DM2}, \cite{AM} and \cite{BNS}.

To the best of our  knowledge, there is no investigation on weakly Laskerian modules over non-Noetherian commutative rings. In
this paper, we investigate some ring-theoretic properties of weakly Laskerian modules over commutative (not necessarily Noetherian)
rings. As a by-product, we deduce several consequences on different types of associated prime ideals. Below, we summarize some of
our main results.

Let $R$ be a weakly Laskerian ring and $I$ an ideal of $R$. We show that:
\begin{enumerate}
\item[i)] ${\rm Min}\ I$ is a finite set; see Theorem \ref{F}.
\item[ii)] If either $\dim R$ is finite or the ring $R[X]$ is weakly Laskerian for some indeterminate $X$ over $R$, then
$\Spec R$ is Noetherain; see Corollary \ref{I} and Theorem \ref{I1}. In particular, in both cases each minimal prime ideal
$\p$ of $I$ is an associated prime of $I$ in the Zariski-Samuel sense; see  Corollary \ref{K}.
\item[iii)] For any weakly Laskerian $R$-module $M$, the trivial ring extension $R\ltimes M$ is weakly Laskerian; see
Theorem \ref{O}.
\item[iv)] The polynomial ring $R[X]$ and the power series ring $R[[X]]$  are not necessarily weakly Laskerian; see Theorem
\ref{T1}.  Thus the analogue of the Hilbert Basis Theorem does not hold for the weakly Laskerianness.
\item[v)] If $A$ is a ring extension of $R$ which is finitely generated as an $R$-module, then $A$ is also a weakly
Laskerian ring; see Theorem \ref{V}.
\end{enumerate}

\section{Minimal prime ideals}

For a proper ideal $I$ of $R$, let ${\rm Min}\  I$ denote the set of all minimal prime ideals of $I$. We know by definition that
if $R$ is a weakly Laskerian ring and $I$ is an ideal of $R$, then the set $\Ass_R R/I$ is finite. But this does not immediately
imply the finiteness of ${\rm Min}\  I$. This is because, it is not true in general that ${\rm Min}\  I\subseteq\Ass_R R/I$.
Let us explain this more.

We start this section by borrowing an example from \cite{An}.

\begin{example}\label{D} Let $R:=\{(a_i)_{i\in \Bbb{N}}\in \underset{i\in \Bbb{N}}\prod \Bbb{Z}/2\Bbb{Z}|\ \ a_i=0\ {\rm for}\,\,
{\rm all\ \ large}\,\, i \ \text{or} \ a_i=1 \ {\rm for}\,\,{\rm all\ \ large}\,\, i \}$. Then with pointwise addition and
multiplication $R$ is a commutative ring with identity. As $x^2=x$ for every $x\in R$, it readily follows that $\Spec R=\Max R$,
and so ${\rm Min }\ (0)=\Spec R$. Let $$\m_i:=\{(a_n)_{n\in
\Bbb{N}}\in R|\  \  a_i=0\}$$ and $$\m_{\infty}:=\{(a_n)_{n\in \Bbb{N}}\in R|\  \ a_n=0\,\, {\rm for}\,\,{\rm all\ \ large}\,\, n\}.$$
It can be easily checked that $\m_{\infty}$ and $\m_i$;  $i\in \mathbb{N}$ are prime and these are the only prime ideals of $R$. Thus
$${\rm Min }\ (0)=\{\m_{\infty}\} \cup(\bigcup_{i\geq 1}\{\m_i\}).$$  For any positive integer $i$, set $\xi_i:=(\delta_{n,i}+2\Bbb{Z})_{n\in
\Bbb{N}}\in R$, where $\delta$ denotes the Kronecker delta. Then, it is easy to see that $\m_i=0:_R\xi_i,$ and hence $\m_i\in
\Ass_R R$ for all positive integers $i$ and that $\m_{\infty}\not\in \Ass_R R$. Thus ${\rm Min }\ (0)\nsubseteq \Ass_RR$. Note
that for any positive integer $i$ we have  $\m_i=(1_{_R}-\xi_i)R$, and so $\m_i$ is finitely generated. Nevertheless, it is easy
to verify that $\m_{\infty}$ is not finitely generated.
\end{example}

In view of the above example, it is  natural to ask: Does any finitely generated minimal prime ideal of $R$ belong to
$\Ass_RR$? The next result gives an affirmative answer to this question.

\begin{proposition}\label{E} Let $I$ be an ideal of $R$. If $\p\in {\rm Min}\ I$ and $\p/I$ is a finitely generated ideal
of the ring $R/I$, then $\p\in \Ass_RR/I$.
\end{proposition}

\proof  Let $\p\in {\rm Min}\ I$ be such that $\p/I$ is a finitely generated ideal of the ring $R/I$. Replacing $R$ with $R/I$,
without loss of generality, we may assume that $I=0$, and so it is enough to show that $\p\in \Ass_RR$.

Since $\p R_{\p}$ is a
finitely generated nilpotent ideal of the ring $R_{\p}$, there exist a positive integer $k$ and an element $s\in R\setminus  \p$
such that $\p^ks=0$. Let $\ell$ be the least positive integer such that $\p^{\ell}s=0$ for some $s\in R\setminus  \p$. So,
$\p^{\ell-1}t\neq 0$ for all $t\in R\setminus  \p$. We claim that $(0:_R\p^{\ell-1}s)=\p$. Assume the contrary. Then, as $\p
\subseteq (0:_R\p^{\ell-1}s)$, there exists an element $s_1\in (0:_R\p^{\ell-1}s) \setminus  \p$. Now, we have $\p^{\ell-1}ss_1=0$
which is a contradiction. So, we have $(0:_R\p^{\ell-1}s)=\p$. Since by the hypothesis $\p$ is finitely generated, it follows that
the ideal $\p^{\ell-1}s$ is also finitely generated, and so $\p=(0:_Ra)$ for some $a\in  \p^{\ell-1}s$. In particular, $\p\in \Ass_RR$.
$\Box$

Concerning Example \ref{D}, we also have the following positive result.

\begin{theorem}\label{F} Let $R$ be a weakly Laskerian ring and $I$ a proper ideal of $R$. Then the set ${\rm Min}\  I$ is finite.
\end{theorem}

\proof As in the proof of Proposition \ref{E}, we may and do assume that $I=0$. So, we should show that the set ${\rm Min }\ (0)$ is finite.

In contrary, assume that ${\rm Min }\ (0)$ is infinite. Then by \cite[Theorem 2.4]{BKN},
there exists an element $\p\in {\rm Min }\ (0)$ such that $\p$ is not finitely generated and for any finitely generated ideal
$J$ of $R$ with $J\subseteq \p$, the set $V(J)\cap {\rm Min }\ (0)$ is infinite.

We inductively choose prime ideals  $\p_1, \p_2,...$ in ${\rm Min }\ (0)\setminus  \{\p\}$ and elements $x_1, x_2, ...$ in $\p$
such that  $x_n\in (\p\cap(\bigcap_{i=1}^{n-1}\p_i))\setminus  \p_n$ and $$\p_n\in \V(Rx_1+Rx_2+\cdots+Rx_{n-1})$$  for all
$n\in \mathbb{N}$. Let $\p_1$ be any element in ${\rm Min }\ (0)\setminus  \{\p\}$ and $x_1$ any element in $\p \setminus  \p_1$.
Next, assume that $n>1$ and prime ideals $\p_1, \p_2,..., \p_{n-1}\in {\rm Min }\ (0)\setminus  \{\p\}$ and elements $x_1, x_2,
...,x_{n-1}\in \p$ with the above requested properties have been chosen. Let $\p_n$ be any element of $$(\V(Rx_1+Rx_2+\cdots
+Rx_{n-1})\cap {\rm Min }\ (0))\setminus  \{\p,\p_1,\p_2,...,\p_{n-1}\}.$$ Then $(\p\cap(\bigcap_{i=1}^{n-1}\p_i))
\nsubseteq \p_n,$ and so we can choose an element $$x_n\in (\p\cap(\bigcap_{i=1}^{n-1}\p_i))\setminus  \p_n.$$ So, the induction
is complete.

For each $n\in \mathbb{N}$, set $I_n:=\p_1 x_1+ \p_2 x_2+\cdots+ \p_n x_n$. Let $n\in \mathbb{N}$ and $0\leq i\leq n$. We show
that $(I_n:x_i)=\p_i$. Clearly, $\p_i\subseteq (I_n:x_i)$. As $\p_i\in \V(Rx_1+Rx_2+\cdots+Rx_{i-1})$ and $x_j\in \p_i$ for all
$j>i$, we deduce that $I_n\subseteq \p_i$. Thus $$(I_n:x_i)x_i\subseteq I_n\subseteq \p_i.$$ But $x_i\notin \p_i$, and so
$(I_n:x_i)\subseteq \p_i$. Set $K:=\bigcup_{n=1}^{\infty}I_n$. Then $$(K:x_i)=\bigcup_{n=1}^{\infty}(I_n:x_i)=\p_i$$
for all $i$. Hence $\p_1,\p_2,...$ are infinitely many associated prime ideals of $R/K$ which is a contradiction.  $\Box$

Recall that a topological space $X$ is said to be Noetherian if any ascending chain of open sets eventually stabilizes.
Refer to \cite[Ch.2, \S 4]{Bo} for more details on Noetherian topological spaces. Our next result provides a criterion for the
Noetherianness of $\Spec R$ equipped with its Zariski topology. We extract the following result from \cite{OP} and apply it several
times in the
sequel.

\begin{lemma}\label{H} The following statements are equivalent:
\begin{enumerate}
\item[i)] $\Spec R$ is Noetherian.
\item[ii)] $\Spec R[X]$ is Noetherian, where $X$ is an indeterminate over $R$.
\item[iii)] $R$ satisfies the ascending chain condition on prime ideals and each ideal has a finite number of minimal
prime ideals.
\item[iv)] Each prime ideal of $R$ is equal to the radical of a finitely generated ideal of $R$.
\end{enumerate}
\end{lemma}

\proof i)$\Leftrightarrow$ii) follows by \cite[Theorem 2.5]{OP} and \cite[Proposition 2.8 (iv)]{OP}.

i)$\Rightarrow$iii) Since $\Spec R$ is Noetherian,  clearly $R$ satisfies the ascending chain condition on prime ideals.
On the other hand, as any closed subset of a Noetherian space has finitely many irreducible components, each ideal of $R$
has a finite number of minimal prime ideals.

iii)$\Rightarrow$i) See \cite[Page 65, Exercise 25]{K}.

i)$\Leftrightarrow$iv) See \cite[Corollary 2.4]{OP}. $\Box$

Our next two results show that, under some mild assumptions, the weakly Laskerianness of $R$ implies the Noetherianness
of $\Spec R$.

\begin{corollary}\label{I} Let $R$ be a finite-dimensional weakly Laskerian ring.  Then $\Spec R$ is Noetherian.
In particular, if $X_1,...,X_n$ are $n$ indeterminacies  over $R$, then $\Spec R[X_1,...,X_n]$ is Noetherian.
\end{corollary}

\proof Since by the hypothesis $R$ has finite dimension, it satisfies the ascending chain condition on prime ideals.
Moreover, in view of Theorem \ref{F} , each ideal of $R$ has a finite number of minimal prime ideals. So by Lemma \ref{H},
$\Spec R$ is Noetherian. The second assertion also follows by Lemma \ref{H}.  $\Box$

\begin{theorem}\label{I1} Let $X$ be an indeterminate over $R$. Assume that the ring $R[X]$ is weakly Laskerian. Then
$\Spec R$ is Noetherian.
\end{theorem}

\proof Suppose the contrary and look for a contradiction. Since the two rings $R$ and $R[X]/XR[X]$ are isomorphic, it follows that
the ring $R$ is also weakly Laskerian. Then, in view of Theorem \ref{F} and Lemma \ref{H}, we deduce that there exists a strictly
increasing chain $$\p_1\subset \p_2\subset \dots \subset\p_n\subset \p_{n+1}\subset \cdots$$ in $\Spec R$. Set $A:=R[X]$ and let
$J$ denote the ideal of $A$ generated by the set $$\{aX^n| \ n\in \mathbb{N} \ \text{and} \ a\in \p_n\}.$$ Also for each natural
integer $n$, set $Q_n:=\p_nA+XA$. Then one may check that $$Q_1\subset Q_2\subset \dots \subset Q_n\subset Q_{n+1}\subset \cdots$$
is a strictly increasing chain in $\Spec A$.

We claim that $\{Q_n\}_{n\in \mathbb{N}}\subseteq \Ass_AA/J$. This will provide the desired contradiction. Let $n\in \mathbb{N}$,
 $b\in \p_{n+1}\setminus \p_n$ and set $c:=bX^n$. We claim that $Q_n=(J:_Ac)$. One has $$Q_nc=b(\p_nX^n)A+(bX^{n+1})A\subseteq J,$$
and so $Q_n\subseteq (J:_Ac)$. Next, let $$h=a_0+a_1X+\cdots +a_tX^t\in (J:_Ac).$$ Then $hbX^n\in J$, and so there are natural
integers $i_1<i_2<\cdots <i_{\ell}$ such that
$$\begin{array}{lll} ba_0X^n+ba_1X^{n+1}+\cdots +ba_tX^{n+t}&=hbX^n\\
&=\sum_{j=1}^{{\ell}}\sum_{k=1}^{n_j}f_{kj}(b_{kj}X^{i_j}),
\end{array}$$
where $b_{kj}\in \p_{i_{j}}$ and $f_{kj}\in A$ for all $j,k$. Comparing the coefficients of $X^n$ in the first and the third terms
of the above display gives $$ba_0\in (b_{kj}|k=1,\dots, n_j, i_j\leq n)R\subseteq \cup_{i=1}^n\p_n=\p_n.$$ As $b\notin \p_n$ one
gets $a_0\in \p_n$, and so $h\in \p_nA+XA=Q_n$. Thus $$Q_n=(J:_Ac)\in \Ass_AA/J. \   \  \Box$$

\begin{definition}(See \cite[Definition 3.1]{HO}.) Let $I$ be an ideal of $R$. A prime ideal $\p$ of $R$ is said to be
an associated prime of $I$ in the Zariski-Samuel sense if $\p=\sqrt{I:_Rx}$ for some $x\in R$. Let ${\rm ZS}(I)$ denote
the set of Zariski-Samuel associated primes of $I$.
\end{definition}

By using Lemma \ref{H}, one can easily deduce the following result:

\begin{lemma}\label{J} Let $I$ be a proper ideal of $R$.  If $\Spec R$ is Noetherian, then ${\rm Min}\ I\subseteq {\rm ZS}(I)$.
\end{lemma}

In view of Lemma \ref{J}, Corollary \ref{I} and Theorem \ref{I1} immediately yield the following result.

\begin{corollary}\label{K} Let $R$ be a weakly Laskerian ring and $I$ a proper ideal of $R$. Assume that either
\begin{enumerate}
\item[i)] $\dim R$ is finite; or
\item[ii)] the ring $R[X]$ is weakly Laskerian for some indeterminate $X$ over $R$.
\end{enumerate}
Then ${\rm Min }\ I\subseteq {\rm ZS}(I)$.
\end{corollary}

Note that by Nagata's celebrated example \cite[Example 1, p 203]{N}, there exist Noetherian integral domains of infinite Krull
dimension. So, the ring $R[X]$ can be weakly Laskerian while $\dim R$ is infinite.

Theorem \ref{F} and Corollary \ref{I} are some instances of the situations in which weakly Laskerian rings behave
like Noetherian rings. However, there are the cases when they behave completely different from Noetherian rings.
See the following example.

\begin{example} Let $R$ be a weakly Laskerian ring and $M$ an $R$-module. One may guess that  $M=0$ if and only if
$\Ass_RM=\emptyset$. Also, one may conjecture that ${\rm Min }\ (0)\subseteq \Ass_RR$. But, the previous two properties
do not hold in general. To this end, let $k$ be a field, $T:=k[X_1,X_2,...]$ and $J:=(X_1,X_2^2,...,X_n^n,...)T$. Let
$R:=T/J$ and $\frak
m=(X_1,X_2,...,X_n,...)R$. Then we have ${\rm Spec}\,R=\{\frak m\},$ and so obviously the ring $R$ is weakly Laskerian.
We claim that ${\rm Ass}_R\, R=\emptyset$. Assume the contrary. Then there is a polynomial $f\in T\setminus J$ such
that $\frak m=0:_R(f+J)$. There exists a positive integer $t$ such that $f\in k[X_1,X_2,...,X_t]$. Then as $f\in T
\setminus J,$ it is easy to see that $(X_{t+1}+J)(f+J)\neq J$ which is a contradiction. Thus $R\neq 0$, $\Ass_RR=
\emptyset$ and  ${\rm Min }\ (0)\nsubseteq \Ass_RR.$
\end{example}

\section{Trivial ring extensions}

Let $M$ be an $R$-module. In this section, we establish a characterization for the weakly Laskerianness of the trivial
ring extension $R\ltimes M$; see Theorem \ref{O}.

Recall that $R\ltimes M:=R\times M$ with addition $(r_1,m_1)+(r_2,m_2):=(r_1 + r_2,m_1 + m_2)$
and multiplication $(r_1,m_1)(r_2,m_2):=(r_1r_2,r_1m_2 + r_2m_1)$ is a commutative ring with identity  $(1,0)$ and is
called the idealization of $M$. Note that $R$ naturally embeds into $R\ltimes M$ via $r \longrightarrow (r,0)$ and if
$N$ is a submodule of $M$, then $0\ltimes N$ is an ideal of $R\ltimes M$. For the ideal $\mathfrak{I}:=0\ltimes M$ of
$R\ltimes M,$  one has $\mathfrak{I}^2=0$. Every ideal of $R\ltimes M$ that contains $0\ltimes M$ has the form $I\ltimes M$
for some ideal $I$ of $R$. In particular, since any prime ideal $\mathfrak{P}$ of $R\ltimes M$ contains all nilpotent
elements of $R\ltimes M$ and hence contains $0\ltimes M$, it follows that $\mathfrak{P}=\p\ltimes M$ for some prime
ideal $\p$ of $R$. Moreover, every ideal of  $R\ltimes M$ that is contained in $0\ltimes M$ has the form $0\ltimes K$
for some submodule $K$ of $M$. Some basic results on idealization can be found in \cite{Hu}.

\cite[Proposition 2.2]{AW} and \cite[Theorem 1.7]{HNN} yield the following characterization for the Noetherianness of
the trivial ring extension $R\ltimes M$.

\begin{proposition}\label{M} Let $M$ be an $R$-module. Then the ring $R\ltimes M$ is Noetherian if and only if the ring $R$
is Noetherian and the $R$-module $M$ is finitely generated.
\end{proposition}

\begin{lemma}\label{N} Let $T$ be a quotient ring of $R$ and $X$ a $T$-module. Then $X$ is weakly Laskerian as an $R$-module
if and only if it is weakly Laskerian as a $T$-module. In particular, $R$ is a weakly Laskerian ring if and only if any quotient
ring of $R$ is  weakly Laskerian.
\end{lemma}

\proof We may assume that $T=R/J$ for some ideal $J$ of $R$. A subset $Y$ of $X$ is a submodule of $X$ as an $R$-module
if and only if it is a submodule of $X$ as a $T$-module. On the other hand, for any $T$-module $Z$ one has $\Ass_RZ \subseteq
V(J)$ and $$\Ass_TZ=\{\frac{\p}{J}| \  \ \p\in \Ass_RZ\}.$$  Thus $|\Ass_RZ|=|\Ass_TZ|$, and so $X$ is weakly
Laskerian as an $R$-module if and only if it is weakly Laskerian as a $T$-module.  $\Box$

\begin{lemma}\label{P} Let $J$ be an ideal of $R$ such that $J^2=0$. Assume that the ring $R/J$ is weakly Laskerian and
the $R/J$-module $J$ is weakly Laskerian. Then the ring $R$ is weakly Laskerian.
\end{lemma}

\proof By Lemma \ref{N}  both $R$-modules $J$ and $R/J$ are weakly Laskerian. Hence by the exact sequence $$0 \rightarrow J
\rightarrow R \rightarrow R/J \rightarrow 0,$$ we deduce that the ring $R$ is weakly Laskerian.  $\Box$

Our next result is the analogues of Proposition \ref{M} for the weakly Laskerianness.

\begin{theorem}\label{O} Let $M$ be an $R$-module. The ring $R\ltimes M$ is weakly Laskerian if and only if $R$ is a weakly
Laskerian ring and $M$ is a weakly Laskerian $R$-module.
\end{theorem}

\proof Set $\mathfrak{J}:=0\ltimes M$. Note that the two rings $R$ and $(R\ltimes M)/\mathfrak{J}$ are naturally isomorphic
and also $\mathfrak{J}$ and $M$ are naturally isomorphic as $R$-modules.

First, assume that $R$ is a weakly Laskerian ring and $M$ is a weakly Laskerian $R$-module. Then Lemma \ref{P} yields that
$R\ltimes M$ is a weakly Laskerian ring.

Conversely, assume that the ring $R\ltimes M$ is weakly Laskerian. As $$R\cong \frac{R\ltimes M}{\mathfrak{J}},$$ it turns out that
the ring $R$ is weakly Laskerian. Moreover, as $\mathfrak{J}$ is a weakly Laskerian $R\ltimes M$-module and $\mathfrak{J}^2=0$,
it follows that $\mathfrak{J}$ is a weakly Laskerian $(R\ltimes M)/\mathfrak{J}$-module. So, $M$ is a weakly Laskerian $R$-module.
$\Box$

As an immediate consequence, we record the following corollary.

\begin{corollary}\label{Q} Let $M$ be a weakly Laskerian module over a Noetherian ring $R$. Then the ring $R\ltimes M$
is weakly Laskerian.
\end{corollary}

We end this section with the following two examples. In the first one, we present a non-Noetherian weakly Laskerian ring.
The second one exhibits a weakly Laskerian ring that is not Laskerian.

In what follows, for an $R$-module $M$, $E_R(M)$ stands for the injective envelope of $M$.

\begin{example}\label{Q1} Let $R$ be a Noetherian semi-local ring and $\p$ a prime ideal of $R$ with $\dim R/\p\leq 1$.
Since $\V(\p)$ is finite, the $R$-module $E_R(R/\p)$ is weakly Laskerian. Hence, Corollary \ref{Q} yields that the ring
$R\ltimes E_R(R/\p)$ is weakly Laskerian. Note that if $\Ht\p>0$, then the $R$-module $E_R(R/\p)$ is not finitely generated,
and so by Proposition \ref{M} the ring $R\ltimes E_R(R/\p)$ is not Noetherian.
\end{example}

\begin{example}\label{Q2} Let $M$ be a Laskerian module and $r$ an element in the Jacobson radical of $R$. Then \cite[Corollary 3.2]{HL}
implies that $\bigcap_{n=1}^{\infty}r^nM=0$. Now, let $(R,\frak m,k)$ be a Noetherian local domain of dimension $d>0$ and let $E:=E_R(k)$.
Then $S:=R\ltimes E$ is a weakly Laskerian local ring with the unique maximal ideal $\frak m \ltimes E$. Let $ 0\neq x\in \frak m$ and
put $r:=(x,0)\in S$. Then $r$ is an element in the Jacobson radical of $S$. Since $xE=E,$  we have $r^nS=x^nR \ltimes E$. In particular,
one has $$0 \ltimes E \subseteq \bigcap_{n=1}^{\infty}r^nS.$$ Thus $\bigcap_{n=1}^{\infty}r^nS\neq 0$, which implies that $S$ is not a
Laskerian ring.
\end{example}

\section{Polynomial ring extensions}

Let $R$ be a weakly Laskerian ring and $\{X_{\gamma}\}_{\gamma\in \Gamma}$ a set of indeterminates over $R$. One may guess
that the rings $R[\{X_{\gamma}\}_{\gamma\in \Gamma}]$ and $R[[\{X_{\gamma}\}_{\gamma\in \Gamma}]]$ are weakly Laskerian.
Theorems \ref {R} and \ref {T1} below show that the finiteness of $\Gamma$ is a necessary but not sufficient condition for
the weakly Laskerianness of these two rings.

\begin{theorem}\label{R} Let $X_{_1},X_{_2},...$ be a countable set of indeterminates over any ring $R$ $($even weakly Laskerian$)$.
Then the rings $R[X_{_1},X_{_2},...]$ and $R[[X_{_1},X_{_2},...]]$ are not weakly Laskerian.
\end{theorem}

\proof We only prove the claim for the ring $R[X_{_1},X_{_2},...]$, because our argument below can be used also for the ring
$R[[X_{_1},X_{_2},...]]$.

Set $A:=R[X_{_1},X_{_2},...]$ and let $\m$ be a maximal ideal of $R$. By Lemma \ref{N}  if $A$ is a weakly Laskerian ring, then
the ring $A/\m A$ is also weakly Laskerian. But, there is
an isomorphism of rings: $$A/\m A\cong (R/\m)[X_{_1},X_{_2},...].$$ As $R/\m$ is a field, it is enough to prove that
the ring $A=k[X_{_1},X_{_2},...]$ is not weakly Laskerian, where $k$ is a field and $X_{_1},X_{_2},...$ are indeterminates
over $k$.

In view of Theorem \ref{F}, it suffices to find an ideal $I$ of $A$ such that the set ${\rm Min}\ I$ is infinite. To this
end, let $$I:=(\{X_{_1}\}\cup\{X_{_2}\}\cup(\bigcup_{n=1}^{\infty}\{X_{_{2^{n}+1}}X_{_{2^{n}+2}}\cdots X_{_{2^{n+1}}} \})).$$
Let $$\mathfrak{B}:=\{\p\in\Spec A|\ \ \p=( X_{_{j_{_1}}}, X_{_{j_{_2}}},X_{_{j_{_3}}},...)\,\,{\rm where}\, j_{_1}=1,j_{_2}=2
\ \ {\rm and} $$ $$\ \  2^{k-2}<j_{_k}\leq 2^{k-1}\ \  {\rm for \ \ all}\ \  k\geq 3\}.$$ Then it is clear that $\mathfrak{B}$
is an infinite subset of $\Spec A$. So, it is enough to prove that ${\rm Min}\ I=\mathfrak{B}$. To do this, first let $\p\in {\rm
Min}\ I$. Then we have $I\subseteq \p$. In particular, $X_{_1},X_{_2}\in \p$ and for each integer $k\geq 3$,
$$X_{_{2^{k-2}+1}}X_{_{2^{k-2}+2}} \cdots X_{_{2^{k-1}}}\in \p$$ which implies $X_{j_{_k}}\in \p$ for some integer
$2^{k-2}<j_{_k}\leq 2^{k-1}$. Now, put $j_{_1}=1$ and $j_{_2}=2$. Then $\q:=(X_{_{j_{_1}}},X_{_{j_{_2}}},X_{_{j_{_3}}},...)$
is a prime ideal of $A$ such that $I\subseteq \q \subseteq \p$. Hence, $\p=\q\in \mathfrak{B}$. Therefore, we have
${\rm Min}\ I\subseteq\mathfrak{B}$. Next, let $\p\in \mathfrak{B}$. Then it is clear that $I\subseteq \p$. So, $\p$
contains a minimal prime ideal $\q$ of $I$. Then $\q \in {\rm Min}\ I\subseteq \mathfrak{B}$. So, $\q\subseteq \p$ and
$\p,\q\in \mathfrak{B}$ which implies $\p=\q$. Note that the elements of $\mathfrak{B}$ are pairwise incomparable under
inclusion.  $\Box$

We will use the following result in the proof of Theorem \ref{T1}. For its proof see \cite[Theorem]{F}.

\begin{proposition}\label{S} Let $X_{_1},X_{_2},...,X_{_n}$ be $n$ indeterminates over $R$. If $A:=R[X_{_1},X_{_2},...,X_{_n}]$,
then $\Ass_AA=\{\p A|\  \ \p\in \Ass_RR \}.$
\end{proposition}

Next, we record the following immediate corollary.

\begin{corollary}\label{T} Let $X_{_1},X_{_2},...,X_{_n}$ be $n$ indeterminates over $R$ and $A:=R[X_{_1},X_{_2},...,X_{_n}]$.
Then for any ideal $I$ of $R$, the two sets $\Ass_RR/I$ and $\Ass_AA/IA$ have the same cardinality. In particular, if the ring
$A$ is weakly Laskerian, then the ring $R$ is weakly Laskerian too.
\end{corollary}

\proof Since $A/IA\cong (R/I)[X_{_1},X_{_2},...,X_{_n}]$, the claim is clear by Proposition \ref{S}.
Note that if $J$ is an ideal of a ring $T$ and $X$ is a $T/J$-module, then $|\Ass_TX|=|\Ass_{T/J}X|$.
$\Box$

\begin{lemma}\label{T2} Let $(R,\m,k)$ be a Noetherian local ring and set $S:=R\ltimes E_R(k)$ and
$B:=S[[X]]$. For any prime ideal $\p$ of $R$, there is a prime ideal $Q\in \Ass_BB$ such that $Q\cap
S=\p\ltimes E_R(k)$.
\end{lemma}

\proof  Let $\p$ be a prime ideal of $R$. As $$E_{R/\p}(k)=(0:_{E_R(k)}\p)=\bigcup_{n=1}^{\infty}(0:_{
E_{R/\p}(k)}\m^n)$$ and for every $n\geq 1$ the $R$-module $(0:_{E_{R/\p}(k)}\m^n)$ is finitely generated,
it follows that the $R$-module $E_{R/\p}(k)$ has a countable generator set $\{a_i\}_{i\in \mathbb{N}_0}$'say.
Now, set $f:=\sum_{i\in
\mathbb{N}_0}(0,a_i)X^i\in B$. As $\Ann_R(E_{R/\p}(k))=\p$, we deduce that the ideal $(0:_{B}f)$ belongs
to the set $$\Omega:=\{J\unlhd B\,|\   \ (0:_{B}f)\subseteq J\,\,{\rm and}\,\,J\cap S=\p\ltimes E_R(k)\}.$$
Because of the natural ring isomorphisms $$\frac{B}{(\p\ltimes E_R(k))[[X]]}\simeq (\frac{S}{\p\ltimes
E_R(k)})[[X]]\simeq (\frac{R}{\p})[[X]],$$ one gets that the ring $B/(\p\ltimes E_R(k))[[X]]$ is Noetherian.
So, it turns out that $\Omega$ has a maximal element $P$. We claim that $P\in \Spec B$. Assume the opposite.
Then there are $\zeta, \xi\in B\backslash P$ such that
$\zeta \xi\in P$. So, by the choose of $P$ there are elements $$x\in (P+B\zeta)\cap S\backslash (\p\ltimes
E_R(k))$$ and $$y\in (P+B\xi)\cap S\backslash (\p\ltimes E_R(k)).$$ Thus $$xy\in (P+B\zeta)(P+B\xi)\cap S
\subseteq P\cap S=\p\ltimes E_R(k)$$ which is a contradiction. So, $P$ is a prime ideal of $B$. Since
$(0:_{B}f)\subseteq
P$, it follows that $P$ contains a minimal prime ideal $Q$ of $(0:_{B}f)$. Now as $$\p\ltimes E_R(k)=
(0:_{B}f)\cap S\subseteq Q\cap S\subseteq P\cap S=\p\ltimes E_R(k),$$ it follows that $Q\cap S=\p\ltimes
E_R(k)$. Since the ring $T:=B/(0:_{B}f)$ is Noetherian and $Q/(0:_{B}f)$ is a minimal prime ideal of $T$,
it follows that $$\frac{Q}{(0:_{B}f)}\in \Ass_TT.$$ Therefore, there is an element $h\in B\backslash (0:_{B}f)$
such that $$Q=((0:_{B}f):_Bh)=(0:_{B}hf),$$ and so $Q\in \Ass_BB$ and $Q\cap S=\p$, as required.  $\Box$

The next result provides an example of a weakly Laskerian ring $S$ such that the rings $S[X]$ and $S[[X]]$
are not weakly Laskerian.

\begin{theorem}\label{T1} Let $(R,\m,k)$ be a Noetherian local ring of dimension $d$ and let $S:=R\ltimes E_R(k)$.
Then the following statements hold:
\begin{enumerate}
\item[i)] $S$ is a weakly Laskerian ring.
\item[ii)] if $d\geq 1$, then the ring $A:=S[X]$ is not weakly Laskerian for any indeterminate $X$ over $S$.
\item[iii)] if $d\geq 2$, then the ring $B:=S[[X]]$ is not weakly Laskerian for any indeterminate $X$ over $S$.
\end{enumerate}
\end{theorem}

\proof i) holds by Example \ref{Q1}.

ii)  As $$\Spec S=\{\p\ltimes E_R(k)\,|\  \ \p \in \Spec R\},$$ it follows that $S$ is a local ring with
the unique maximal ideal $\n:=\m\ltimes E_R(k)$. In addition for the ideal $\mathfrak{J}:=0\ltimes E_R(k)$ of $S$, we
have $\mathfrak{J}^2=0$, and so for the ideal $J:=\mathfrak{J}[X]$ of the polynomials ring $A:=S[X]$ we have $J^2=0$.
So, $J$ has an $A/J$-module structure. But, by the ring isomorphisms $$\frac{A}{J}\simeq (\frac{S}{\mathfrak{J}})[X]
\simeq R[X],$$  it turns out that $A/J$ is a Noetherian ring. Now, we claim that the ring $A$ is not weakly Laskerian.
In contrary assume that the ring $A$ is weakly Laskerian. Then the ideal $J$ of $A$ is a weakly Laskerian $A$-module
and hence by the $A/J$-module structure of $J$, it follows that $J$ is also a weakly Laskerian $A/J$-module. Hence,
by \cite[Theorem 3.3]{Ba}, the $A/J$-module $J$ is an  ${\rm FSF}$ module. So, by the definition there exists a finitely
generated submodule $L$ of $J$ such that the $A/J$-module $J/L$ has finite support. But, in this situation $L$ is a finitely
generated ideal of $A$. Then there are elements $f_1, f_2,..., f_n\in J$ such that $L=(f_1, f_2,..., f_n)A.$ Next, let
$f_i=\Sigma_{j=0}^{k_i}(0,b_{i j})X^j$ for $i=1,2,..., n$. Then $B:=\Sigma_{i=1}^n\Sigma_{j=0}^{k_i}Rb_{i j}$ is a finitely
generated submodule of the Artinian $R$-module $E:=E_R(k)$. Since $d\geq 1$, the $R$-module $E$ is not finitely generated, and
hence $\Ass_RE/B=\{\m\}$. Moreover it is obvious that $L\subseteq (0\ltimes B)[X]$, and so
the $A/J$-module $\frac{J}{(0\ltimes B)[X]}$ has finite support. Thus, the $A$-module $\frac{J}{(0\ltimes B)[X]}$ has finite
support too. Since $\m \in \Ass_RE/B$ it easily follows that $\n \in \Ass_S(\frac{\mathfrak{J}}{(0\ltimes B)})$, and so by
Proposition \ref {S}, we have $\n[X]\in \Ass_A(\frac{J}{(0\ltimes B)[X]})$. This implies that $$V(\n[X])\subseteq \Supp_A(
\frac{J}{(0\ltimes B)[X]}).$$
Since the PID $k[X]$ has infinitely many non-associated irreducible elements, it becomes clear that $\Spec k[X]$ is
infinite. Hence, from the natural ring isomorphisms $$\frac{A}{\n[X]}\simeq (\frac{S}{\n})[X] \simeq k[X],$$ we deduce
that $V(\n[X])$ is an infinite subset of $\Spec A$. Now, we have achieved the desired contradiction.

iii) Since $\dim R=d\geq 2$, it follows that $\Spec R$ and consequently $\Spec S$ is finite. By Lemma \ref{T2}, for any
prime ideal $\p$ of $R$, there is a prime ideal $Q\in \Ass_BB$ such that $Q\cap S=\p\ltimes E_R(k)$. Thus the finiteness
of $\Ass_BB$ implies the finiteness of $\Spec S$. Therefore, the ring $B$ is not weakly Laskerian. $\Box$

\section{Integral ring extensions}

Theorem \ref{V} below is the main result of this section. To prove it, we need the following result which might be of independent
interest.

\begin{lemma}\label{U} Let $X$ be an indeterminant over $R$ and $A:=R[X]$. Let $J$ be an ideal of $A$, $\q\in \Ass_AA/J$ and
$\p=\q\cap R$.
For each integer $k\geq 0$, let $\ab_k$ denote the set of all $a\in R$ for which there exists a polynomial of the type $a_0+a_1X+\cdots
+a_{k-1}X^{k-1}+aX^k$ in $J$. Then $\ab_0\subseteq \ab_1\subseteq \ab_2\subseteq \cdots $ is a chain of ideals of $R$ and there
exists an integer $n\geq 0$ such that $\p\in \Ass_RR/\ab_n$.
\end{lemma}

\proof It is easy to see that  $\ab_0\subseteq \ab_1\subseteq \ab_2\subseteq \cdots $ is a chain of ideals of $R$.
By the definition, there is an element $f \in A\setminus  J$ such that $\q=(J:_Af)$. We can choose an element $f=a_0
+a_1X+\cdots+a_nX^n\in A$ of the minimum degree with the property $\q=(J:_Af)$. Next, we claim that $(\ab_n:_Ra_n)=\p$.
Assume the contrary. Then as $\p \subseteq (\ab_n:_Ra_n)$,  there is an element $a\in (\ab_n:_Ra_n)\setminus \p$.
As $a\in R$ and $a\not\in \q\cap R=\p$, it follows that $a\not\in \q$.  Since $aa_n\in \ab_n$, there exists $g\in J$
of degree at most $n$ such that the degree of $af-g$ is less than $n$. As $$\q= (J:_Af)\subseteq (J:_Aaf)=J:_A(af-g),$$
by the choose of $f$, it follows that $\q\subsetneqq (J:_Aaf)$. Hence, there exists an element $h\in (J:_Aaf)\setminus
\q$. Now, we have $ha\in (J:_Af)=\q$, $h\not\in \q$ and $a\not\in \q$ which is a contradiction. Thus we have
$(\ab_n:_Ra_n)=\p$, and so $\p\in \Ass_RR/\ab_n$. $\Box$

\begin{theorem}\label{V} Let $R$ be a weakly Laskerian ring and $A$ a ring extension of $R$ which is finitely generated
as an $R$-module. Then $A$ is also a weakly Laskerian ring.
\end{theorem}

\proof There are elements $\theta_1,...,\theta_n\in A$ such that $A=R[\theta_1,...,\theta_n]$ and $\theta_i$'s are integral over
$R$. As $$R[\theta_1,\theta_2, ..., \theta_n]=(R[\theta_1,\theta_2, ..., \theta_{n-1}])[\theta_n],$$ by induction on $n$, we may
assume that $A=R[\theta]$ and $\theta$ is integral over $R$.

Let $X$ be an indeterminant over $R$ and define $\phi:R[X]\rightarrow A$ with $$\phi(c_0+c_1X+\cdots+c_kX^k)=c_0+c_1\theta+
\cdots+c_k\theta^k.$$ Then $\phi$ is a surjective ring homomorphism, and so $A\cong R[X]/J$, where $J:={\rm ker}(\phi)$.
Hence, it is enough to prove that $T:=R[X]/J$ is a weakly Laskerian ring. To this end, let $J_1/J$ be an ideal of $T$. We
have to show that the set $\Ass_{T}(R[X]/J_1)$ is finite. Set ${\frak a}:=J_1\cap R$. Then by \cite[Proposition 5.6 i)]{AM},
the extension $R/{\frak a}\subseteq R[X]/J_1$ is integral and finitely generated. Since $\theta$ is integral over $R$,
there exists a polynomial $$a_0+a_1X+\cdots+a_{t-1}X^{t-1}+X^t\in J\subseteq J_1.$$ For each integer $k\geq 0$, set $$\ab_k:
=\{a\in R|\ \ {\rm there}\,\,{\rm exists}\,\,{\rm an}\,\,{\rm element}\,\,a_0
+a_1X+\cdots+a_{k-1}X^{k-1}+aX^k\in J_1\}.$$ Then, $\ab_0\subseteq \ab_1\subseteq \ab_2\subseteq \cdots $ is a chain of
ideals of $R$, and as $1_{_R}\in \ab_t$ it follows that $R=\ab_t=\ab_{t+1}=\ab_{t+2}=\cdots$.

Assume that $\Ass_{R[X]}(\frac{R[X]}{J_1})$ is infinite. Set $$\mathfrak{D}:=\{\q\cap R|\ \q\in \Ass_{R[X]}(\frac{R[X]}{J_1})\}.$$
Then as $R=\ab_t=\ab_{t+1}=\ab_{t+2}=\cdots$, Lemma \ref{U} implies that $\mathfrak{D}\subseteq\bigcup_{n=0}^{t-1}\Ass_RR/\ab_n$.
In particular, as $R$ is a weakly Laskerian ring, we deduce that $\mathfrak{D}$ is a finite set. So, there exists
an element $\p\in \mathfrak{D}$ such that there are infinitely many elements in $\Ass_{R[X]}(\frac{R[X]}{J_1})$
lying over $\p$. But for each $\q\in \Ass_{R[X]}(\frac{R[X]}{J_1})$ with $\q\cap R=\p$, $\q/J_1$ is a prime ideal
of $R[X]/J_1$ lying over $\p/{\frak a}$. So, there are infinitely many prime ideals of $R[X]/J_1$ lying over
$\p/{\frak a}$. But this a contradiction with \cite[Exercise 9.3]{Mat}. So, $A$ is a weakly Laskerian ring. $\Box$

As an easy conclusion, we bring  the following result.

\begin{corollary}\label{Z4}  Let $R$ be a weakly Laskerian ring and $X$ an indeterminant over $R$. Let $J$ be an ideal of
the ring $R[X]$ which contains a monic polynomial $f$. Then the ring $R[X]/J$ is weakly Laskerian.
\end{corollary}

\proof The ring $R[X]/J$ is a quotient of the ring $A:=R[X]/fR[X]$. As $A$ is a finitely generated ring extension of $R$, the
claim follows by Theorem \ref{V}.  $\Box$

We end the paper with the following result.

\begin{proposition}\label{W} Let $R$ be a Noetherian ring and $M$ a weakly Laskerian $R$-module. Let $T$ be a Noetherian
semi-local $R$-algebra which is integral over $R$. Then $M\otimes_{R}T$ is a weakly Laskerian $T$-module.
\end{proposition}

\proof By \cite[Theorem 3.3]{Ba}, there exists a finitely generated submodule $N$ of $M$ such that ${\rm {\rm Supp}}_{_{R}}\,M/N$
is finite, and so ${\rm dim}_{_{R}}M/N\leq 1$. Set $J:=\bigcap_{\p\in {\rm {\rm Supp_R}}\,M/N}\p$. Then ${\rm dim}\, R/J\leq 1$
and ${\rm {\rm Supp}}_R\,M/N=V(J)$. It is easy to check that $${\rm Supp}_{T}\,(M/N\otimes_RT)\subseteq V(JT).$$ Because $T$ is
integral over $R$, \cite[Proposition 5.6 i)]{AM} yields that $T/JT$ is integral over $R/JT\cap R$, and so we deduce that
$${\rm dim}\,\frac{T}{JT}={\rm dim}\,\frac{R}{JT\cap R}\leq {\rm dim}\,\frac{R}{J}\leq 1.$$ So as $T$ is a semi-local ring, $V(JT)$
and, consequently, ${\rm {\rm {\rm Supp}}}_{T}\,(M/N\otimes_RT)$ is finite. Now by applying \cite[Theorem 3.3]{Ba} again,  we
conclude that $M\otimes_RT$ is a weakly Laskerian $T$-module.
$\Box$

\begin{acknowledgement}  We thank the referee for his/her valuable comments and suggestions on the paper.
\end{acknowledgement}


\end{document}